# On the benefits of output sparsity for multi-label classification

Evgenii Chzhen*[1,2,3], Christophe Denis†[1], Mohamed Hebiri‡[1], and Joseph Salmon§[2]

[1]LAMA, Université Paris-Est Marne-la-Vallée, Université Paris-Est, 5 Boulevard Descartes, 77420, Champs-sur-Marne, France
[2]LTCI, Télécom ParisTech, Université Paris-Saclay, 46 rue Barrault, 75013, Paris, France,
[3]Institute for Information Transmission Problems of the Russian Academy of Sciences (Kharkevich Institute), Bolshoy Karetny per. 19, build.1, 127051, Moscow, Russia

**Abstract**

The multi-label classification framework, where each observation can be associated with a set of labels, has generated a tremendous amount of attention over recent years. The modern multi-label problems are typically large-scale in terms of number of observations, features and labels, and the amount of labels can even be comparable with the amount of observations. In this context, different remedies have been proposed to overcome the curse of dimensionality. In this work, we aim at exploiting the output sparsity by introducing a new loss, called the sparse weighted Hamming loss. This proposed loss can be seen as a weighted version of classical ones, where active and inactive labels are weighted separately. Leveraging the influence of sparsity in the loss function, we provide improved generalization bounds for the empirical risk minimizer, a suitable property for large-scale problems. For this new loss, we derive rates of convergence linear in the underlying output-sparsity rather than linear in the number of labels. In practice, minimizing the associated risk can be performed efficiently by using convex surrogates and modern convex optimization algorithms. We provide experiments on various real-world datasets demonstrating the pertinence of our approach when compared to non-weighted techniques.

## 1 Introduction

The field of Extreme Multi-Label Classification (EMLC), is concerned with classification problems where the number of classes is large and possibly comparable with the number of observations, with one specific characteristic: a data point is tagged with several labels. It has recently attracted a vast amount of contributions due to the variety of problems that EMLC can model: text categorization [Gao et al., 2004], functional genomics [Barutcuoglu et al., 2006], image classification [Li et al., 2014] to name a few. The objective in multi-label classification is to predict a binary vector $Y \in \{0,1\}^L$ for a given observation $X \in \mathcal{X}$, where $L$ is the number of labels available and $\mathcal{X} = \mathbb{R}^D$ is the feature space.

A pioneering approach, called Binary Relevance (BR), requires to build $L$ binary classifiers [Tsoumakas et al., 2010] in a one-versus-all fashion. Yet, when the number of classes is large, possibly of the same order as the available observations, BR becomes too expensive. Moreover, BR degrades when the labels are highly correlated and poor performance were established [Dembczyński et al., 2012] in this context. Such limitations led to different workarounds such as probabilistic classifier chains [Cheng et al., 2010], graph based algorithms [Tan et al., 2015] and correlated logistic regression [Li et al., 2016].

*evgenii.chzhen@univ-paris-est.fr
†christophe.denis@u-pem.fr
‡mohamed.hebiri@u-pem.fr
§joseph.salmon@telecom-paristech.fr

Our motivation in this work starts from the observation that for a given $X \in \mathcal{X}$ the output label vector $Y$ is often sparse: for each data point only a few labels are generally active in real-world scenarios (*i.e.,* $Y$ has few non-zero coordinates, say $K$, with $K \ll L$). To the best of our knowledge, exploiting the output sparsity for improved learning appeared first in [Hsu et al., 2009] for EMLC, where the authors relied on compressed sensing [Donoho, 2006, Candès et al., 2006] to approximate the output by $K \log(D)$ measurements. In terms of theory Hsu et al. [2009] defined the notion of sparsity for $\mathbb{E}[Y|X]$. Hence the observed outcomes might not necessary be sparse. Meanwhile, in several real-world applications of EMCL the output vectors are indeed sparse, see Table 3. Of course this can be due to errors in the labeling process: some meaningful labels are expected to be missed by human annotators when the number of labels $L$ is huge. Consequences of such errors were investigated in [Jain et al., 2016], where the sparsity was defined in terms of relevant labels for a each data point.

Here, we aim at providing, in a principle way, a statistical framework for EMLC by exploiting ideas of sparsity and some heuristics proposed in [Hsu et al., 2009, Jain et al., 2016]. We introduce a novel sparsity assumption, which reflects the nature of multi-label classification problems. Our assumption, unlike the one proposed in [Hsu et al., 2009], restricts the possibility to observe a large amount of ones in the output labels for any observed data point. We also propose a new loss, tailored to handle output sparsity. The proposed loss function aims at producing a sparse prediction with the possibility to make a small amount of mistakes on the zero labels of the underlying signal. Such a feature is also addressed by Jain et al. [2016] who proposed to omit zero labels to build their classifier. This allows to handle scenarios where the annotator is not malicious, and does not add wrong tags: this is typically the case for a dataset like Wikipedia where the community would automatically erase wrong labels, while missing ones might be common due to the millions of possible tags. Another benefit of our approach is that the naive classifier that outputs only ones, is not achieving the best performance, as it is the case in [Jain et al., 2016].

The organization of the paper is as follows. In Section 2, we introduce our framework, our output sparsity model and our new loss function for EMLC. We also discuss the advantages of the proposed loss over standard ones. In Section 3, we adapt classical generalization bounds for this loss. These bounds are advantageous due to their linear dependence on the output sparsity level $K$. As a consequence of our analysis we introduce a novel notion of Rademacher complexity. It can be seen as the classical Rademacher complexity separated into two parts, associated with either active labels (the ones) or inactive labels (the zeros). Our analysis is somehow similar to the one introduced in [Yu et al., 2014], though our assumption and the loss function of Section 2 allows us to develop faster rates of convergence. The final contribution of this section is an adaptation of Zhang's lemma [Zhang, 2004] to this context. In Section 4.1, we perform a numerical study in order to show the advantage of our loss function over classical ones, for both synthetic and real data.

## 2 Framework and Notation

Throughout the paper, the following notation is used: for any $N, L, D \in \mathbb{N}$, $[D] := \{1, \dots, D\}$, we write $N \vee L := \max\{N, L\}$ and $N \wedge L := \min\{N, L\}$. We use $N$ for the number of observations, $L$ for the number of labels and $D$ for the dimension of the feature space $\mathcal{X}$. For a matrix $\boldsymbol{Q}$, we denote by $\boldsymbol{Q}^\top$ its transpose. The Euclidean norm is written $\|\cdot\|_2$ and for any $a, b \in \mathbb{R}^M$ we denote by $\langle a, b\rangle = a^\top b$ their inner product. For any $\boldsymbol{A}, \boldsymbol{B} \in \mathbb{R}^{N\times L}$ we denote by $\langle \boldsymbol{A}, \boldsymbol{B}\rangle = \mathrm{Tr}(\boldsymbol{A}^\top \boldsymbol{B})$ their Hilbert-Schmidt inner product. For a given matrix $\boldsymbol{A} \in \mathbb{R}^{N\times L}$ and for any $p \geq 1$, its Schatten $p$-norm [Bhatia, 1997] is written $\|\boldsymbol{A}\|_{\sigma,p} := \left(\sum_{i=1}^{L\wedge N} \sigma_i^p(\boldsymbol{A})\right)^{1/p}$, with $\sigma_i(\boldsymbol{A})$ being the singular values of $\boldsymbol{A}$ in decreasing order. We use $\|\boldsymbol{A}\|_{\sigma,\infty} := \sigma_1(\boldsymbol{A})$ to denote the spectral norm of the matrix $\boldsymbol{A}$. We write $X \in \mathcal{X} = \mathbb{R}^D$ for a feature vector and $Y \in \{0,1\}^L$ for its associated label vector ($Y^l$ being the $l^{\text{th}}$ component of $Y$), and $\mathbb{P}$ be the underlying probability distribution of $\mathcal{X} \times \{0,1\}^L$.

**Definition 1** (Sparse vectors)**.** We say that a binary vector $Y \in \{0,1\}^L$ is $K$-sparse for $K \in [L]$ if

$$\|Y\|_0 := \#\{l \in [L] : Y^l = 1\} = K \quad , \tag{1}$$

where $\#A$ refers to the cardinality of a finite set $A$.

Sparsity is an important property that often arises in statistics and signal processing to bypass the curse of dimensionality, especially in high dimensional regression or inverse problems [Elad, 2010, Bühlmann and van de Geer, 2011]. In this work, though, we are rather interested in sparsity of the output variable $Y$, a somehow unconventional context.

**Assumption 1** (Label sparsity)**.** There exists a positive integer $K \ll L$, *s.t.* ,

$$\mathbb{P}\left\{\sum_{l=1}^{L} \mathbb{1}_{\{Y^l=1\}} > K | X\right\} = 0, \text{ a.s.} \quad .$$

We emphasize that $K \ll L$ and restrict our attention only to this extreme case. Note that Assumption 1 is different from the one presented in [Hsu et al., 2009] and is in fact stronger. Indeed, under Assumption 1, $\mathbb{E}[\sum_{l=1}^{L} \mathbb{1}_{\{Y^l=1\}} | X]$ is bounded by $K$. Meanwhile, assuming that $\mathbb{E}[\sum_{l=1}^{L} \mathbb{1}_{\{Y^l=1\}} | X] \leq K$ does not necessary mean that every realization of $Y$ would be sparse.

**Definition 2** (Active labels)**.** We say that a label $l \in [L]$ is active for a binary vector $Y \in \{0,1\}^L$ if $Y^l = 1$. If $Y^l = 0$ we say that $l$ is inactive. The set of all active (inactive) labels is called active (inactive) set of the vector $Y$.

The choice of the loss function is crucial for EMLC. One of the most popular one is the Hamming loss, however this loss does not take into account label sparsity, since it treats all labels equally. In [Jain et al., 2016], the authors proposed propensity loss functions, that exploit the observation that the label vector $Y$ for a given data point $X$, often does not consists of all possible active labels. Being more precise, the authors introduced another label vector $Y^*$, which is a noisy version of the true label vector $Y$. Assuming, that the true label vector $Y$ is deterministic, the authors defined a label to be relevant, if it is active in the underlying deterministic label vector $Y$.

Here we consider the following observation framework. Let $\mathcal{D}_N = \{(X_1, Y_1), \ldots, (X_N, Y_N)\}$ be *i.i.d.* copies of $(X, Y)$. Let $\mathbb{P}_N$ be the empirical probability distribution. Let $\mathcal{F}$ be a class of hypothesis, *i.e.*, each element $f \in \mathcal{F}$ is a measurable function from $\mathcal{X}$ to $\mathbb{R}^L$. With weights $p_0, p_1 \geq 0$ such that $p_0 + p_1 = 1$ we can consider the following asymmetric risk

$$\mathcal{R}(f) = \sum_{l=1}^{L} \left\{ p_0 \mathbb{E}\left[\mathcal{L}_0(f^l(X)) \mathbb{1}_{\{Y^l=0\}}\right] + p_1 \mathbb{E}\left[\mathcal{L}_1(f^l(X)) \mathbb{1}_{\{Y^l=1\}}\right] \right\} \quad . \tag{2}$$

The empirical counterpart of the risk in Eq. (2) is given by

$$\hat{\mathcal{R}}_N(f) = \frac{1}{N} \sum_{i=1}^{N} \sum_{l=1}^{L} \left\{ p_0 \mathcal{L}_0(f^l(X_i)) \mathbb{1}_{\{Y_i^l=0\}} + p_1 \mathcal{L}_1(f^l(X_i)) \mathbb{1}_{\{Y_i^l=1\}} \right\} \quad . \tag{3}$$

Note that we also write $\hat{\mathcal{R}}_N(f, \mathcal{D}_N)$ for $\hat{\mathcal{R}}_N(f)$ when the dependency on $\mathcal{D}_N$ is needed. In this work we are interested in the empirical risk minimizer (ERM)

$$\hat{f} \in \arg\min_{f \in \mathcal{F}} \hat{\mathcal{R}}_N(f) \quad , \tag{4}$$

| | $\hat{Y}_0 \equiv 0$ | $\hat{Y}_1 \equiv 1$ | $\hat{Y}_{2K}$ | $\hat{Y}_w \equiv 1 - Y$ |
|---|---|---|---|---|
| Proposed loss | $p_1 K$ | $p_0(L-K)$ | $p_0 K$ | $p_1 K + p_0(L-K)$ |
| Hamming loss: $p_0 = p_1 = 0.5$ | $\frac{K}{2}$ | $\frac{L-K}{2}$ | $\frac{K}{2}$ | $\frac{L}{2}$ |
| [Jain et al., 2016]: $p_0 = 0, \quad p_1 = 1$ | $K$ | $0$ | $0$ | $K$ |
| Our choice: $p_0 = \frac{2K}{L}, \quad p_1 = 1 - \frac{2K}{L}$ | $K - \frac{2K^2}{L}$ | $2K - \frac{2K^2}{L}$ | $\frac{2K^2}{L}$ | $3K - \frac{4K^2}{L}$ |

Table 1: Loss examples costs for several classifiers with underlying true label being $K$-sparse, with $K \ll L$: $\hat{Y}_0 \equiv 0$: output no label, $\hat{Y}_1 \equiv 1$ output all labels, $\hat{Y}_{2K}$: output correct active set plus $K$ mistakes on inactive set, $\hat{Y}_w \equiv 1 - Y$: always wrong

and to compare $\hat{f}$ with $f^*$, the best predictor[1] in $\mathcal{F}$

$$f^* \in \arg\min_{f \in \mathcal{F}} \mathcal{R}(f) \ . \tag{5}$$

## 2.1 Motivation and examples

In this section we provide some motivational examples for the risk in Equation (2). First, we notice that this risk is essentially asymmetric, and is motivated by the fact that under Assumption 1 the output label vector consists of at most $K$ active labels. Treating $Y^l = 1$ and $Y^l = 0$ in the same way might force classifiers to overfit on inactive labels only.

We handle this phenomenon by adding weights to our loss function, and such weights only rely on the (global) sparsity of the labels. This approach is different from the one presented in [Jain et al., 2016], where the authors proposed to add weights that are (local) label dependent. Intuitively $p_0$ should be a small weight to reflect the presence of inactive labels and $p_1$ should be almost 1 to give priority to active labels over inactive ones.

A particular case of our loss includes $\mathcal{L}_j(\hat{Y}^l) = \mathcal{L}_j^{0/1}(\hat{Y}^l) = \mathbb{1}_{\{Y^l \neq j\}}$ for $j \in \{0, 1\}$, where $\hat{Y}$ is a binary prediction vector. This leads to introducing a weighted Hamming loss written as

$$\mathcal{L}_{0/1}^{w}(Y, \hat{Y}) = \sum_{l=1}^{L} \left\{ p_0 \mathbb{1}_{\{\hat{Y}^l = 1\}} \mathbb{1}_{\{Y^l = 0\}} + p_1 \mathbb{1}_{\{\hat{Y}^l = 0\}} \mathbb{1}_{\{Y^l = 1\}} \right\} \ . \tag{6}$$

We remind that the Hamming loss can be written as

$$\mathcal{L}^{H}(Y, \hat{Y}) = \sum_{l=1}^{L} \mathbb{1}_{\{Y^l \neq \hat{Y}^l\}} \ . \tag{7}$$

Let us now study the advantage of the weighted Hamming loss with specific weights.

One of the drawbacks of the Hamming loss is that it treats all labels equally and therefore, does not take into account label sparsity. For more insight, let us consider the scenario where $Y$ is exactly $K$-sparse and let us analyze the following classifiers (a synthesis is also given in Table 1):

- $\hat{Y}_0 \equiv 0$: predicts all labels inactive,
- $\hat{Y}_1 \equiv 1$: predicts all labels active,
- $\hat{Y}_w \equiv 1 - Y$: msspredicts all labels,
- $\hat{Y}_{2K}$: correctly predicts the active set of $Y$ and makes exactly $K$ mistakes on its inactive set.

Intuitively, one would like to build a loss which is able to differentiate between the first three predictions and $\hat{Y}_{2K}$. Indeed, in large-scale problems the predictions similar to $\hat{Y}_{2K}$ provide more valuable insights compared to the first three ones.

[1] One does not have access to $f^*$ since $\mathbb{P}$ is unknown.

With our choice of weights

$$p_0 = \frac{2K}{L}, \quad p_1 = 1 - \frac{2K}{L} \quad , \tag{8}$$

the introduced loss function treats $\hat{Y}_0, \hat{Y}_1$ and $\hat{Y}_w$ almost equally and it promotes predictions with small amount of mistakes on inactive sets and correct predictions on active sets. Meanwhile, the Hamming loss does not make any difference between $\hat{Y}_{2K}$ and $\hat{Y}_0$, and the loss considered in [Jain et al., 2016] gives a high promotion to naive classifiers like $\hat{Y}_1$.

**Remark 1.** In practice, the weights in Eq. (8) rely on the unknown sparsity constant $K$ in Assumption 1. Since this quantity is unknown to the practitioner, a simple strategy consists in performing a rough estimation based on the observed labels. Hence, we consider

$$\hat{p}_0 = \frac{2\hat{K}}{L}, \quad \hat{p}_1 = 1 - \hat{p}_0 \quad , \tag{9}$$

where we estimate the output sparsity level by the maximal sparsity on the observations:

$$\hat{K} = \max_{i \in [N]} \sum_{l=1}^{L} Y_i^l \quad .$$

## 3 Generalization Bounds

In this section we study generalization properties of empirical risk minimization (ERM) under sparsity Assumption 1. The generalization bounds usually have the following form: with high probability, the error of the estimator is bounded by an empirical estimate of the error plus a residual term depending on the complexity of the hypothesis class investigated [Bartlett et al., 2005].

There are various notions of complexities, one of the most popular being the Vapnik-Chervonenkis (VC) dimension [Vapnik, 1998]. However, VC-dimension is distribution-free, hence model independent, so it is not well tailored to our context. An important alternative technique to estimate complexity of the hypothesis class relies on Rademacher averages, which provide data driven bounds. In the context of multilabel classification it was first used by [Yu et al., 2014] that additionally accounted for possible missing labels.

In Section 3.1 we focus on adapting standard analysis based on Rademacher complexities [Koltchinskii, 2006, Bartlett et al., 2002, Bartlett and Mendelson, 2002] to provide generalization bounds linear in $K$. In Section 3.1.1 we adapt Talagrand's concentration inequality in the form of [Bousquet, 2002]. This result is a first step towards bounds based on local Rademacher complexities [Bartlett et al., 2005], but further development on this road is outside the scope of this paper. In Section 3.2 we prove a generalization of Zhang's Lemma [Zhang, 2004] for the case of weighted loss functions.

### 3.1 Rademacher complexities

The standard technique of comparing the ERM with the best classifier in the class $\mathcal{F}$ [Bartlett and Mendelson, 2002] relies on the following inequality

$$\mathcal{R}(\hat{f}) \leq \mathcal{R}(f^*) + 2 \sup_{f \in \mathcal{F}} \left| \mathcal{R}(f) - \hat{\mathcal{R}}_N(f) \right| \quad . \tag{10}$$

From now on, w.l.o.g. we assume $\mathcal{L}_0, \mathcal{L}_1$ are bounded by 1. Hence, our next step is to notice that the function $\sup_{f \in \mathcal{F}} \left| \mathcal{R}(f) - \hat{\mathcal{R}}_N(f, \mathcal{D}_N) \right|$ has bounded differences and that McDiarmid's inequality [McDiarmid, 1989], reminded in the supplementary material, can be applied to control the supremum in terms of its expectation. With Assumption 1, we can state

**Proposition 1.** *Consider the following function*

$$\varphi(\mathcal{D}_N) = \sup_{f\in\mathcal{F}} \left|\mathcal{R}(f) - \hat{\mathcal{R}}_N(f,\mathcal{D}_N)\right| \ ,$$

*and let $\mathcal{D}'_N$ be a second set of observations such that $\mathcal{D}'_N$ and $\mathcal{D}_N$ differ only in one observation $(X_k, Y_k)$ and $(X'_k, Y'_k)$, where $(X'_k, Y'_k)$ is a copy of $(X, Y)$ independent from $\mathcal{D}_N$. Then,*

$$\left|\varphi(\mathcal{D}_N) - \varphi(\mathcal{D}'_N)\right| \le \frac{L}{N} p_0 + \frac{2K}{N} p_1, \ a.s. \ . \tag{11}$$

*Proof.* See supplementary material for the proof. □

Notice that choosing $p_0 = 0, p_1 = 1$ leads to the best upper bound in terms of $K$. However, this choice of weights is useless in practice: $\hat{Y}_1$ the naive classifier that always predicts a label as active, is a minimizer of the empirical risk. To avoid this degeneracy, we recommend $p_0 > 0$. In particular, the proposed choice in Equation (8) allows to bound the differences of $\varphi$ almost surely by a factor $4K/N$.

Without Assumption 1 the differences would be bounded by a constant of order $L/N$, a possibly not informative bound in large-scale scenarios with $N$ and $L$ being large.

We also emphasize that McDiarmid's inequality consists in bounding differences of the supremum over the whole domain. However in our case, we can get a better bound almost surely, following [Combes, 2015]. In the light of the above discussion we have the following lemma

**Lemma 1.** Setting weights as in Eq. (8) and $\delta \in ]0,1[$, then with probability at least $1-\delta$ the following inequality holds

$$\sup_{f\in\mathcal{F}} \left|\mathcal{R}(f) - \hat{\mathcal{R}}_N(f)\right| \le \mathbb{E}_N \sup_{f\in\mathcal{F}} \left|\mathcal{R}(f) - \hat{\mathcal{R}}_N(f)\right| + 4\frac{K}{\sqrt{N}}\sqrt{\log\frac{1}{\delta}} \ .$$

Using a standard symmetrization inequality to control the expected deviation of the empirical mean from the true one in terms of Rademacher averages, we can get the following proposition

**Proposition 2.** *For any class of hypothesis $\mathcal{F}$ of measurable functions*

$$\mathbb{E}_N \sup_{f\in\mathcal{F}} \left|\mathcal{R}(f) - \hat{\mathcal{R}}_N(f)\right| \le 2p_0\mathfrak{R}_0(\mathcal{L}_0 \circ \mathcal{F}) + 2p_1\mathfrak{R}_1(\mathcal{L}_1 \circ \mathcal{F}) \ ,$$

*where for $j = 0$ and $j = 1$*

$$\mathfrak{R}_j(\mathcal{F}) = \frac{1}{N}\mathbb{E}_N\mathbb{E}_\varepsilon \sup_{f\in\mathcal{F}} \left|\sum_{i=1}^{N}\sum_{l=1}^{L} \varepsilon_i^l f^l(X_i)\mathbb{1}_{\{Y_i^l = j\}}\right| \ , \tag{12}$$

*are two Rademacher complexities of the class $\mathcal{F}$ and $\varepsilon_i^l$ are* i.i.d. *Rademacher variables.*

Assuming that $\mathcal{L}_0$ and $\mathcal{L}_1$ are $C$-Lipschitz and using contraction property [Ledoux and Talagrand, 1991, Theorem 4.12.] we can write

$$\mathbb{E}_N \sup_{f\in\mathcal{F}} \left|\mathcal{R}(f) - \hat{\mathcal{R}}_N(f)\right| \le 2Cp_0\mathfrak{R}_0(\mathcal{F}) + 2Cp_1\mathfrak{R}_1(\mathcal{F}) \ .$$

Combining Lemma 1, Proposition 2 and the contraction property [Ledoux and Talagrand, 1991, Theorem 4.12.] we get the following theorem

**Theorem 1.** *Setting weights as in Eq. (8) and $\delta \in ]0,1[$, and assuming that $\mathcal{L}_0$ and $\mathcal{L}_1$ are $C$-Lipschitz, then with probability at least $1-\delta$, the following bound holds*

$$\mathcal{R}(\hat{f}) \leq \mathcal{R}(f^*) + 4Cp_0\mathfrak{R}_0(\mathcal{F}) + 4Cp_1\mathfrak{R}_1(\mathcal{F}) + 8\frac{K}{\sqrt{N}}\sqrt{\log\frac{1}{\delta}} \ .$$

Notice that, the previous theorem provides rates with linear dependence on $K$, hence our next goal is to control $\mathfrak{R}_1(\mathcal{F})$ and $\mathfrak{R}_0(\mathcal{F})$, with possibly the following rate w.r.t. to $K, L$ and $N$

$$\mathfrak{R}_0(\mathcal{F}) \lesssim \frac{L}{\sqrt{N}}, \quad \mathfrak{R}_1(\mathcal{F}) \lesssim \frac{K}{\sqrt{N}} \ .$$

Such bounds are needed in order to have linear dependency on $K$. Indeed this would lead to $p_1\mathfrak{R}_1(\mathcal{F}) \simeq \mathfrak{R}_1(\mathcal{F}) \lesssim K/\sqrt{N}$ and $p_0\mathfrak{R}_0(\mathcal{F}) \simeq K\mathfrak{R}_0(\mathcal{F})/L \lesssim K/\sqrt{N}$.

#### 3.1.1 Refined bound: Talagrand's inequality

Similar but somehow sharper bounds on the risk of the empirical minimizer $\hat{f}$ can be obtained if we manage to bound the largest variance of the class members $f \in \mathcal{F}$ [Bartlett et al., 2005]. This technique relies on Talagrand's concentration inequality [Bousquet, 2002] for empirical processes. The only difference from the previous analysis holds in a refined version of Lemma 1.

**Lemma 2.** Setting weights as in Eq. (8), $\delta \in ]0,1[$ and letting $r > 0$ be such that

$$\forall j \in \{0,1\}, \quad \forall f \in \mathcal{F}, \quad \mathbb{E}[\max_{l\in[L]} \mathcal{L}_j^2(f^l(X))] \leq r \ , \tag{13}$$

then with probability at least $1-\delta$ the following inequality holds

$$\sup_{f\in\mathcal{F}} \left|\mathcal{R}(f) - \hat{\mathcal{R}}_N(f)\right| \leq 2\mathbb{E}\sup_{f\in\mathcal{F}} \left|\mathcal{R}(f) - \hat{\mathcal{R}}_N(f)\right| + K\sqrt{\frac{32r}{N}\log\frac{1}{\delta} + \frac{10K}{3N}\log\frac{1}{\delta}} \ .$$

*Proof.* See supplementary material for details. □

Applying Proposition 2 and the contraction principle [Ledoux and Talagrand, 1991, Theorem 4.12.] we get the following theorem

**Theorem 2.** *Setting weights as in Eq. (8), $\delta \in ]0,1[$ and letting $r > 0$ be such that the Conditions (13) are satisfied, we have with probability at least $1-\delta$*

$$\mathcal{R}(\hat{f}) \leq \mathcal{R}(f^*) + 8Cp_0\mathfrak{R}_0(\mathcal{F}) + 8Cp_1\mathfrak{R}_1(\mathcal{F}) + K\sqrt{\frac{128r}{N}\log\frac{1}{\delta} + \frac{40K}{3N}\log\frac{1}{\delta}} \ .$$

It is important to point out that this bound does not give an improvement in $K$, however it can be made tighter if the variance bound $r$ decreases w.r.t. $N$. As mentioned in [Bartlett et al., 2005], Theorem 2 is not useful when applied to the whole class of functions $\mathcal{F}$. Meanwhile, by applying Theorem 2 to a subset $\mathcal{F}' \subset \mathcal{F}$ or to a modified version of $\mathcal{F}$, a better result could be obtained [Bartlett et al., 2005, Theorem 3.3.]. In multi-label settings, this type of bounds was first considered in [Xu et al., 2016], where the authors developed a specific regularization for ERM-type algorithms. Though, they used local Rademacher complexities and did not investigate the influence of $K$ on their bound.

### 3.2 Margin-based performance bound

Previous sections have focused on establishing bounds on the risk (2) of some estimator $\hat{f}$ that lies in some class of functions $\mathcal{F}$. Here, we are more concerned by controlling the (weighted) excess risk of classifiers, a more natural goal in multi-label setting.

A classifier is a measurable function $g : \mathcal{X} \to \{0,1\}^L$, and here we restrict our attention to classifiers of the form

$$g_f^l(x) = \begin{cases} 1, & f^l(x) \geq 0 \\ 0, & \text{otherwise} \end{cases} , \tag{14}$$

where[2] $f : \mathcal{X} \mapsto \mathbb{R}^L$. A natural measure of performance in the context of multi-label classification with sparse labels is the Weighted Hamming risk discussed in Section 2.1

$$\mathcal{R}_{0/1}^w(f) = \sum_{l=1}^{L} \left\{ p_0 \mathbb{E}\left[\mathbb{1}_{\{f^l(X)\geq 0\}}\mathbb{1}_{\{Y^l=0\}}\right] + p_1 \mathbb{E}\left[\mathbb{1}_{\{f^l(X)<0\}}\mathbb{1}_{\{Y^l=1\}}\right]\right\} , \tag{15}$$

with $\mathcal{L}_0^{0/1}(f^l) = \mathbb{1}_{\{f^l(X)\geq 0\}}$, $\mathcal{L}_1^{0/1}(f^l) = \mathbb{1}_{\{f^l(X)<0\}}$ and weights $p_0, p_1$ as in Eq. (8). Let $\eta(X) = (\eta^1(X), \ldots, \eta^L(X))^\top$ be the regression function, where each component is given by $\eta^l(X) = \mathbb{P}\{Y^l = 1|X\}$. The Bayes classifier is given by

**Proposition 3.** *The minimizer of Eq. (15) is given for all $l \in [L]$ by*

$$(g_\eta^*)^l(x) = \begin{cases} 1, & \eta^l(x) \geq p_0 \\ 0, & \eta^l(x) < p_0 \end{cases} .$$

We denote by $\mathcal{R}_{0/1}^{w*} = \mathcal{R}_{0/1}^w(\eta^l(x) - p_0)$ the minimum of Eq. (15) and define the Weighted Hamming excess risk of a classifier $f$ as

$$\mathcal{E}_{0/1}^w(f) = \mathcal{R}_{0/1}^w(f) - \mathcal{R}_{0/1}^{w*} . \tag{16}$$

Since the cost functions $\mathcal{L}_j^{0/1}(f^l)$ used in the risk (15) are not $C$-Lipschitz, the analysis provided in the previous sections is not applicable. The purpose of the present section is to develop a different analysis to bound the above excess risk.

Given the set of observations $\mathcal{D}_N$, an estimator of the Bayes classifier can be defined as the minimizer of the empirical counterpart of the above risk

$$\hat{\mathcal{R}}_{N,0/1}^w(f) = \frac{1}{N}\sum_{i=1}^{N}\sum_{l=1}^{L}\left\{p_0 \mathbb{1}_{\{f^l(X_i)\geq 0\}}\mathbb{1}_{\{Y^l=0\}} + p_1 \mathbb{1}_{\{f^l(X_i)<0\}}\mathbb{1}_{\{Y^l=1\}}\right\} . \tag{17}$$

Minimization the former is a computationally difficult problem, since indicator functions are not convex. Therefore some convex surrogate of $\mathcal{L}_j^{0/1}(f^l)$ can be introduced. Establishing a generalization of Zhang's Lemma [Zhang, 2004] for the case of weighted loss functions, we are able to show good performance in terms of the Weighted Hamming excess risk Eq. (16) while dealing with its convexification [Zhang, 2004, Bartlett et al., 2006]. Due to major advances in convex optimization theory [Boyd and Vandenberghe, 2004] it has become possible to solve large-scale convex problems efficiently. Convex optimization in machine learning and more specifically in the classification literature has proved its efficiency. It is reflected by the popularity of such methods as Boosting [Freund and Schapire, 1997], logistic regression [Friedman et al., 2000] and support vector machine [Vapnik, 1998].

[2] for instance, $f$ can be defined as in the previous sections

Let $\phi : \mathbb{R} \mapsto \mathbb{R}_+$ be a non-negative convex surrogate. We consider a convexified form of Equation (15)

$$\mathcal{R}_\phi^w(f) = \sum_{l=1}^{L} \left\{ p_0 \mathbb{E}\Big[\phi(-f^l(X))\mathbb{1}_{\{Y^l=0\}}\Big] + p_1 \mathbb{E}\Big[\phi(f^l(X))\mathbb{1}_{\{Y^l=1\}}\Big] \right\} \ , \tag{18}$$

where $f : \mathcal{X} \mapsto \mathbb{R}^L$ (one may keep in mind to Equation (14) to deduce a classifier $g$ from $f$). Now, the aim is to show that optimizing an empirical version of a convexified risk can be beneficial in terms of weighted Hamming loss. For convenience, we also introduce the following notation $\mathcal{R}_\phi^w(\eta^l, f^l) = p_1\eta^l\phi(f^l) + p_0(1-\eta^l)\phi(-f^l)$, so that

$$\mathcal{R}_\phi^w(f) = \sum_{l=1}^{L} \mathbb{E}_X\Big[\mathcal{R}_\phi^w(\eta^l(X), f^l(X))\Big] \ .$$

We write $f_\phi^*(\eta) : \mathcal{X} \mapsto \mathbb{R}^L$ for the minimizer

$$f_\phi^*(\eta) \in \arg\min_f \mathcal{R}_\phi^w(\eta, f) \ , \tag{19}$$

and $\Delta\mathcal{R}_\phi^w(\eta^l, f^l)$ for the label-wise excess risk

$$\Delta\mathcal{R}_\phi^w(\eta^l, f^l) = \mathcal{R}_\phi^w(\eta^l, f^l) - \mathcal{R}_\phi^w(\eta^l, (f_\phi^*)^l) \ .$$

Now, we are ready to state the main result of this section

**Lemma 3** (Weighted Zhang's Lemma)**.** Assume that for all $l \in [L]$ the function $\phi$ is classification-calibrated, *i.e.*, $(f_\phi^*)^l(\eta^l) > 0$ when $\eta^l > p_0$ and $(f_\phi^*)^l(\eta^l) < 0$ when $\eta^l < p_0$. Assume there exist $c > 0$ and $s \geq 1$ such that for all $l \in [L]$ and for all $\eta^l \in [0,1]$,

$$\left|p_0 - \eta^l\right|^s \leq c^s \Delta\mathcal{R}_\phi^w(\eta^l, 0) \ . \tag{20}$$

Then for any measurable vector function $f$,

$$\mathcal{R}^w(f) - \mathcal{R}^{w*} \leq cL^{\frac{s-1}{s}} \Big(\mathcal{R}_\phi^w(f) - \mathcal{R}_\phi^w(f_\phi^*)\Big)^{\frac{1}{s}} \ . \tag{21}$$

*Proof.* See supplementary material for details. □

Notice that in the context of multi-label classification, convexification deteriorates the rate of convergence in terms of its dependence on the number of labels. While in previous sections the rate was linear in $K$, we now pay an additional $L^{(s-1)/s}$ term (replacing the total number of labels $L$ by the sparsity constant $K$ is a point we plan to address in future works). The term on the right hand side $\mathcal{R}_\phi^w(f) - \mathcal{R}_\phi^w(f_\phi^*)$ can be decomposed into estimation and approximation errors, both depending on the complexity of the class $\mathcal{F}$. However, the approximation error could also be bounded at the rate $K/N$ up to a log factor (depending on the complexity of $\mathcal{F}$) by adapting more refined techniques, see for instance in [Bartlett et al., 2006, Theorem 4]. Since our loss function is asymmetric we assume $(f_\phi^*)^l(\eta^l) < 0$ when $\eta^l < p_0$ in Lemma 3. In the original paper [Zhang, 2004] the author assumed that the surrogate function $\phi$ was invariant with respect to reflection of coordinates. Note that standard surrogate functions satisfy all the assumptions of Lemma 3:

**Example 1.** Here we provide some examples of convex functions $\phi$ which satisfy the assumptions of Lemma 3:

| Settings $N = 2D = 2L = 200$ | Median output sparsity | | Recall (micro) | | Precision (micro) | |
|---|---|---|---|---|---|---|
| | Weighted | Classical | Weighted | Classical | Weighted | Classical |
| $K = 2$ | 2.47(0.02) | 0.04(0.0) | 1.0(0.0) | 0.02(0.0) | 0.80(0.0) | 1.0(0.03) |
| $K = 6$ | 6.83(0.01) | 0.43(0.01) | 1.0(0.0) | 0.07(0.0) | 0.88(0.0) | 1.0(0.0) |
| $K = 10$ | 9.85(0.03) | 1.81(0.04) | 0.90(0.0) | 0.18(0.0) | 0.91(0.0) | 1.0(0.0) |
| $K = 14$ | 10.90(0.14) | 4.11(0.15) | 0.72(0.0) | 0.29(0.0) | 0.93(0.0) | 0.99(0.0) |
| $K = 18$ | 10.98(0.21) | 6.61(0.2) | 0.58(0.0) | 0.36(0.0) | 0.95(0.0) | 0.99(0.0) |

Table 2: Numerical experiments with synthetic data. Evaluation of each quantity is averages over 100 simulations and the variance is reported in brackets. The regularization parameter is fixed $\lambda = 30$.

- Square loss : $\phi(v) = (1-v)^2$ with $s = 2, c = \sqrt{2}$ and

$$(f^*_\phi)^l(\eta^l) = \frac{\eta^l - p_0}{\eta^l(1-2p_0) + p_0},$$

- Boosting loss: $\phi(v) = e^{-v}$ with $s = 2, c = \sqrt{2}$ and

$$(f^*_\phi)^l(\eta^l) = \frac{1}{2}\log\frac{(1-p_0)\eta^l}{(1-\eta^l)p_0}.$$

*Proof.* See supplementary material for details. □

# 4 Linear prediction approach and experiments

We now provide explicit bounds on Rademacher complexities for linear approximation of $Y$. This setting has been recently studied by Yu et al. [2014], Xu et al. [2016], Chen and Lin [2012], though the choice of loss function considered by the authors does not handle output sparsity. We start with a generalization bound for the linear prediction case.

**Theorem 3.** *Let $\mathcal{X} = \mathbb{R}^D$ and $f(X_i) = \boldsymbol{W}^\top X_i$, where $\boldsymbol{W} \in \mathbb{R}^{D\times L}$, $\mathcal{W} = \{\boldsymbol{W} \in \mathbb{R}^{D\times L} : \|\boldsymbol{W}\|_{\sigma,1} \leq \lambda\}$. Let $\boldsymbol{W}^*$ be the minimum of the risk Eq. (2) over $\mathcal{W}$, and let $\hat{\boldsymbol{W}}$ be the minimum of its empirical counterpart Eq. (3). Assume that $\mathbb{E}_X[\|X\|_2] \leq 1$, and weights as in Eq. (8). Then, with probability at least $1-\delta$ we have*

$$\mathcal{R}(\hat{\boldsymbol{W}}) \leq \mathcal{R}(\boldsymbol{W}^*) + 4\lambda C\left(\sqrt{\frac{K}{N}} + \frac{K}{\sqrt{NL}}\right) + 8\frac{K}{\sqrt{N}}\sqrt{\log\frac{1}{\delta}}\ .$$

*Proof.* See supplementary material for details. □

Theorem 3 is an important case of the analysis developed in Section 3.1, illustrating that with convex surrogates, we can still get rate of convergence linear in $K$. Also notice that the term $K/\sqrt{NL}$ is small as compared to $\sqrt{K/N}$, due to our sparse scenario $K \ll L$.

This result also provides a theoretical way to choose the parameter $\lambda \simeq \sqrt{K}$ controlling the size of $\mathcal{W}$. By increasing $\lambda$ and, as a consequence the set $\mathcal{W}$, we obtain a better approximation of the Bayes rule in Proposition 3. At the same time the increase of $\lambda$ leads to a worse estimation error in Theorem 3.

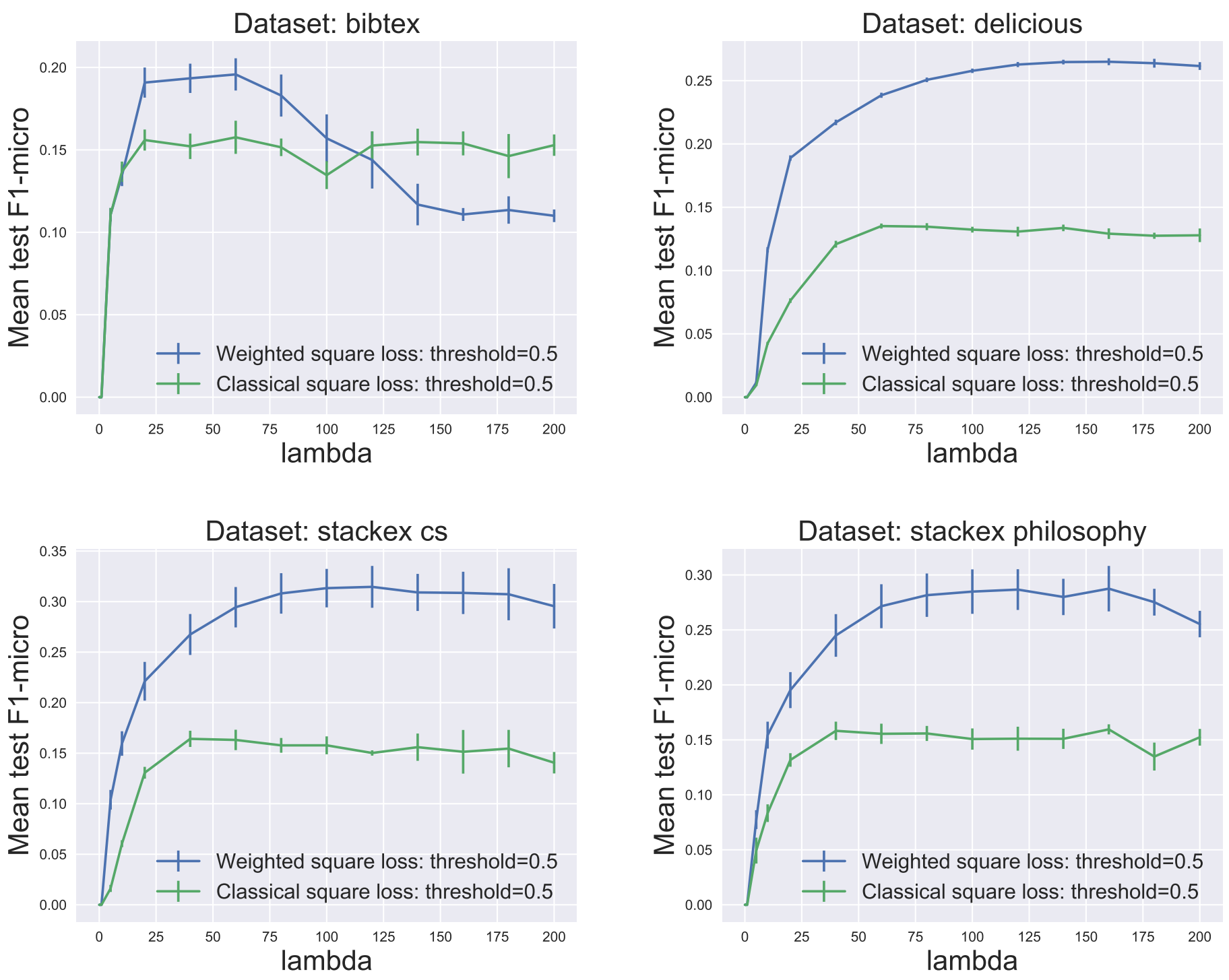

Figure 1: Mean test $F_1$-micro score for various datasets, with $\theta$ fixed to 0.5 for predicting with Equation (23).

## 4.1 Experiments

Let us introduce some notation to compute the estimator analyzed in Theorem 3. First, let $\boldsymbol{W} = [W^1, \ldots, W^L] \in \mathbb{R}^{D \times L}$, $\mathcal{W} = \{\boldsymbol{W} \in \mathbb{R}^{D \times L} : \|\boldsymbol{W}\|_{\sigma,1} \leq \lambda\}$ and $X_i \in \mathbb{R}^D$ for each $i \in [N]$. We denote by $\boldsymbol{X} = [X_1, \ldots, X_N]^\top \in \mathbb{R}^{N \times D}$ the design matrix and by $\boldsymbol{Y} = (Y_i^l)_{i=1,\ldots,N, l=1,\ldots,L} \in \mathbb{R}^{N \times L}$ the label matrix. The estimator we considered is the solution of the convex problem

$$\hat{\boldsymbol{W}} \in \underset{\boldsymbol{W} \in \mathcal{W}}{\arg\min} \frac{1}{N} \left\| \hat{\boldsymbol{P}} \odot (\boldsymbol{Y} - \boldsymbol{X}\boldsymbol{W}) \right\|_{\sigma,2}^2 , \tag{22}$$

where we denote by $\boldsymbol{A} \odot \boldsymbol{B}$ the Hadamard product of two matrices, and by $\hat{\boldsymbol{P}} \in \mathbb{R}^{N \times L}$ the weight matrix defined by

$$\hat{P}_i^l = \begin{cases} \sqrt{\hat{p}_0} & \text{if } Y_i^l = 0 \\ \sqrt{\hat{p}_1} & \text{if } Y_i^l = 1 \end{cases} ,$$

Notice that by choosing $\hat{p}_0 = \hat{p}_1 = 1/2$, $\hat{\boldsymbol{W}}$ becomes a usual Frobenius norm minimizer under nuclear norm constraint [Xu et al., 2016, Yu et al., 2014, Chen and Lin, 2012]. However our choice of the weight matrix $\hat{\boldsymbol{P}}$ allows to handle sparse label space. The minimization problem in (22) is solved with a Frank-Wolf algorithm [Jaggi and Sulovský, 2010, Jaggi, 2013], reminded in supplementary material, since it is a simple and efficient algorithm, which requires to compute only the largest singular value of the gradient at each iteration. Once the matrix $\hat{\boldsymbol{W}}$ is obtained, we form the final prediction by thresholding each entry of the matrix $\boldsymbol{X}\hat{\boldsymbol{W}}$, *i.e.,* for a threshold $\theta$

$$\hat{\boldsymbol{Y}}_i^l = \mathbb{1}_{\{\langle X_i, \hat{W}^l \rangle \geq \theta\}} . \tag{23}$$

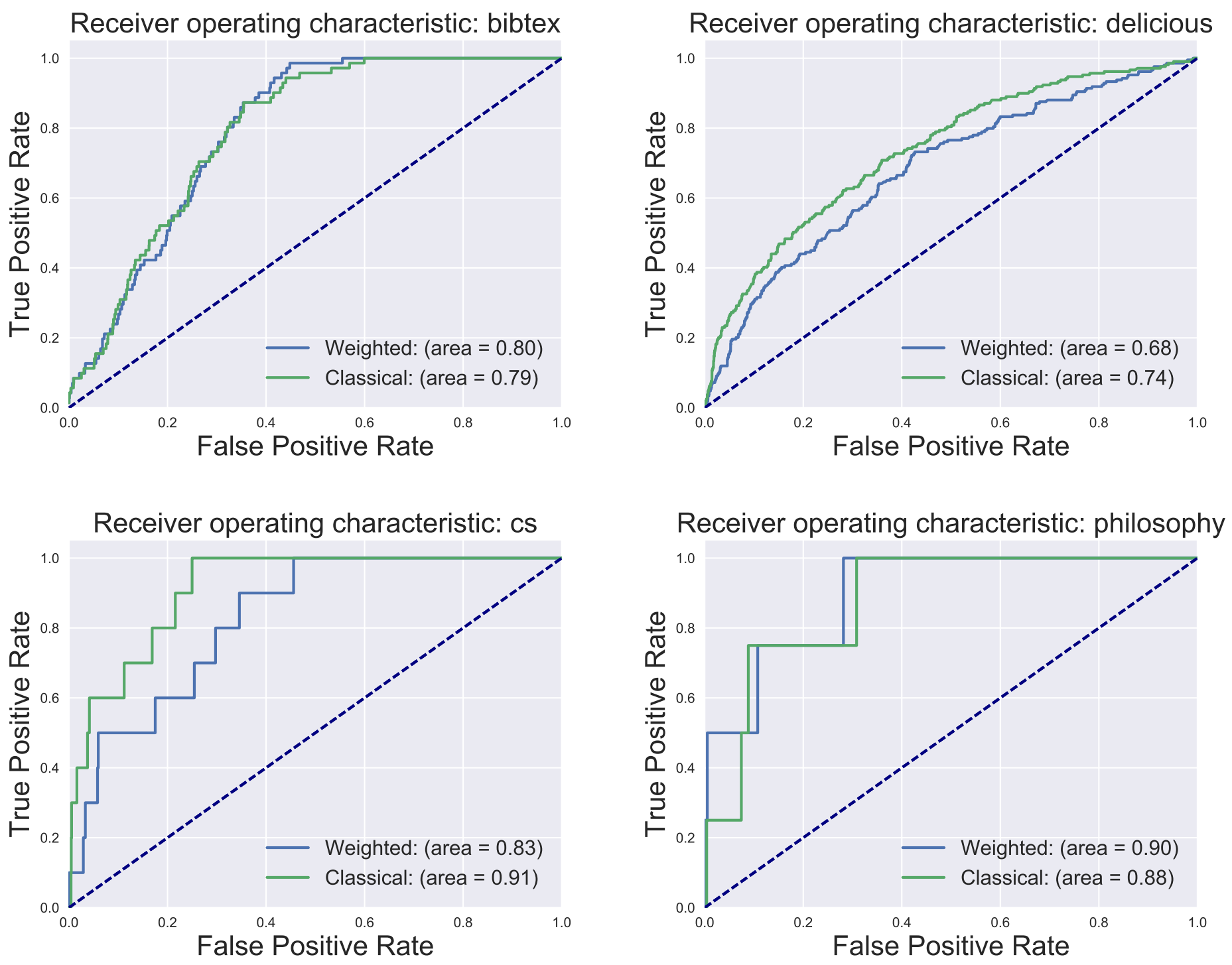

Figure 2: Micro averaged ROC curve for the best classifier obtained via cross-validation over $\lambda$ in Equation (22).

To get a relevant regularization parameter, we have performed 5-fold cross-validation on each real-world dataset, and for each scoring function with regularizers selected among the following values: $\lambda \in \{0.1, 0.5, 1, 10, 20, 40, 60, 80, 100, 120, 140, 160, 180, 200\}$. To evaluate the influence of the loss functions on the method, we have considered the following score functions, see for instance [Natarajan and Jain, 2016]:

- $F_1$-micro score,
- AUC (micro averaged),
- Recall (micro averaged),
- Precision (micro averaged).

The goal of our experiments is to show the advantage of the weighted scheme proposed in Equation (22) over non-weighted techniques, in terms of different performance measures. The benefit of our introduced loss is reflected in

## 4.2 Data

For our empirical study we consider both synthetic and real-world datasets. By controlling the ration $K/L$ on synthetic dataset (see supplementary material for details), we can show empirical evidence to support our framework. We also considered the following real-world datasets, *bibtex* [Katakis et al., 2008], *delicious* [Tsoumakas et al., 2008], *stackex_cs, _philosophy* [Charte et al., 2015]. The real-world datasets are obtained via *mldr* [Charte and Charte, 2015] R package, and normalized to satisfy the assumption of the Theorem 3. Table 3 presents a description of real-world datasets with a focus on the estimated sparsity.

| Dataset | $N$ | $D$ | $L$ | $\hat{K}$ | $\hat{p}_0$ |
|---|---|---|---|---|---|
| bibtex | 7395 | 1836 | 159 | 28 | 0.352 |
| delicious | 16105 | 500 | 983 | 25 | 0.051 |
| stackex_philosphy | 3971 | 842 | 233 | 5 | 0.043 |
| stackex_cs | 9270 | 635 | 274 | 5 | 0.036 |

Table 3: Datasets summary: Number of observations ($N$), features dimensions ($D$), number of labels ($L$), maximal observed sparsity ($\hat{K}$) and $\hat{p}_0$ evaluated as Eq. (9).

### 4.3 Results

The performance of the proposed method on synthetic data is reported in Table 2. We observe a decreasing performance of our method as $K$ grows, but note that for the large sparsity constants (*e.g.,* $K \geq 14$), our scenario $K \ll L$ no longer holds.

The performance of out algorithm on the real-world datasets is reported on Figure 1 and 2. The proposed algorithm achieves better performance over classical (non-weighted) in terms of $F_1$-micro score. It is important to point out that our introduced loss is admitting a prediction to be good if it has few amount of mistakes on the inactive set of the underlying signal $Y$. Since, ROC curves consider only True Positive Rates and False Positive Rates, the method we proposed may show worse performance in terms of AUC. Indeed, by allowing some amount of mistakes on the inactive set of the vector $Y$, we increase the False Positive Rate and hence, decrease the performance in terms of the AUC score. This effect is illustrated on Figure 2. To plot the ROC curves, we consider the classifier with $\lambda$ fixed achieving the best performance in terms of $F_1$-micro score, and each point on the curve corresponds to different threshold levels $\theta$ in Eq. (9).

## 5 Conclusion

We have introduced a new framework for multi-label classification with large number of classes. We have defined a natural output sparsity assumption and provided an associated weighted loss (and risk) leveraging this sparsity property. Interestingly, we could prove generalization bounds linear in the sparsity of the labels. We could also control the performance of a minimizer of a convexified version of our proposed risk. This leads to good practical performance, not only for the introduced risk but also for standard scores in multi-label setting such as recall or $F_1$-scores.

## 6 Acknowledgement

We would like to thank Jean-Baptiste Schiratti, Mónika Csikós and Guillaume Papa for their useful comments and remarks. This work was partially supported by ”laboratoire d'excellence Bézout of Université Paris Est” and by ”Chair Machine Learning for Big Data at Télécom ParisTech”.

## References

P. L. Bartlett and S. Mendelson. Rademacher and Gaussian complexities: risk bounds and structural results. *J. Mach. Learn. Res.*, 3:463–482, 2002. 5

P. L. Bartlett, S. Boucheron, and G. Lugosi. Model selection and error estimation. *Mach. Learn.*, 48(1-3): 85–113, 2002. 5

P. L. Bartlett, O. Bousquet, and S. Mendelson. Local Rademacher complexities. *Ann. Statist.*, 33(4): 1497–1537, 2005. 5, 7

P. L. Bartlett, M. I. Jordan, and J. D. Mcauliffe. Convexity, classification, and risk bounds. *J. Am. Statist. Assoc.*, 101(473):138–156, 2006. 8, 9

Z. Barutcuoglu, R. E. Schapire., and O. G. Troyanskaya. Hierarchical multi-label prediction of gene function. *Bioinformatics*, 22(7):830–836, 2006. 1

R. Bhatia. *Matrix analysis*, volume 169 of *Graduate Texts in Mathematics*. Springer-Verlag, New York, 1997. 2

O. Bousquet. A Bennett concentration inequality and its application to suprema of empirical processes. *C. R. Acad. Sci. Paris, Ser. I*, 334:495–500, 2002. 5, 7, 16

S. Boyd and L. Vandenberghe. *Convex optimization*. Cambridge University Press, Cambridge, 2004. 8

P. Bühlmann and S. van de Geer. *Statistics for high-dimensional data*. Springer Series in Statistics. Springer, Heidelberg, 2011. Methods, theory and applications. 3

E. J. Candès, J. Romberg, and T. Tao. Robust uncertainty principles: Exact signal reconstruction from highly incomplete frequency information. *IEEE Trans. Inf. Theory*, 52(2):489–509, 2006. 2

F. Charte and D. Charte. Working with multilabel datasets in R: The mldr package. *The R Journal*, 7(2): 149–162, dec 2015. 12

F. Charte, A. J. Rivera, M. J. del Jesus, and F. Herrera. QUINTA: a question tagging assistant to improve the answering ratio in electronic forums. In *EUROCON*, pages 1–6, 2015. 12

Y.-N. Chen and H.-T. Lin. Feature-aware label space dimension reduction for multi-label classification. In *NIPS*, pages 1529–1537, 2012. 10, 11

W. Cheng, E. Hüllermeier, and K. Dembczynski. Bayes optimal multilabel classification via probabilistic classifier chains. In *ICML*, pages 279–286, 2010. 2

R. Combes. An extension of McDiarmid's inequality. *arXiv preprint arXiv:1511.05240*, 2015. 6

K. Dembczyński, W. Waegeman, W. Cheng, and E. Hüllermeier. On label dependence and loss minimization in multi-label classification. *Mach. Learn.*, 88(1-2):5–45, 2012. 1

D. L. Donoho. Compressed sensing. *IEEE Trans. Inf. Theory*, 52(4):1289–1306, 2006. 2

M. Elad. *Sparse and redundant representations. From theory to applications in signal and image processing*. New York, NY: Springer, 2010. 3

M. Frank and P. Wolfe. An algorithm for quadratic programming. *Naval research logistics quarterly*, 3 (1-2):95–110, 1956. 21

Y. Freund and R. E. Schapire. A decision-theoretic generalization of on-line learning and an application to boosting. *Journal of computer and system sciences*, 55(1):119–139, 1997. 8

J. Friedman, T. Hastie, and R. R. Tibshirani. Additive logistic regression: a statistical view of boosting. *Ann. Statist.*, 28(2):337–407, 2000. 8

S. Gao, W. Wu, C. H. Lee, and T. S. Chua. A MFoM learning approach to robust multiclass multi-label text categorization. In *ICML*, pages 329–336, 2004. 1

D.J. Hsu, M. Sham Kakade, J. Langford, and T. Zhang. Multi-label prediction via compressed sensing. In *NIPS*, pages 772–780, 2009. 2, 3

M. Jaggi. Revisiting Frank-Wolfe: Projection-free sparse convex optimization. In *ICML*, pages 427–435, 2013. 11, 21, 22

M. Jaggi and M. Sulovský. A Simple Algorithm for Nuclear Norm Regularized Problems. *ICML*, 2010. 11, 21

H. Jain, Y. Prabhu, and M. Varma. Extreme multi-label loss functions for recommendation, tagging, ranking & other missing label applications. In *KDD*, pages 935–944, 2016. 2, 3, 4, 5

I. Katakis, G. Tsoumakas, and I. Vlahavas. Multilabel text classification for automated tag suggestion. In *Proc. ECML/PKDD 2008 Discovery Challenge*, pages 75–83, 2008. 12

V. Koltchinskii. Rademacher penalties and structural risk minimization. *IEEE Trans. Inf. Theory*, 47(5): 1902–1914, September 2006. 5

M. Ledoux and M. Talagrand. *Probability in Banach spaces. Isoperimetry and processes.* Springer-Verlag, Berlin, 1991. 6, 7

Q. Li, B. Xie, J. You, W. Bian, and D. Tao. Correlated logistic model with elastic net regularization for multilabel image classification. *IEEE Trans. Image Process.*, 25(8):3801–3813, Aug 2016. 2

X. Li, F. Zhao, and Y. Guo. Multi-label image classification with a probabilistic label enhancement model. In *UAI*, pages 430–439, 2014. 1

C. McDiarmid. On the method of bounded differences. In *Surveys in Combinatorics*, number 141 in London Mathematical Society Lecture Note Series, pages 148–188. Cambridge University Press, 1989. 5, 16

N. Natarajan and P. Jain. Regret bounds for non-decomposable metrics with missing labels. In *NIPS*, pages 2874–2882. 2016. 12

M. Tan, Q. Shi, A. van den Hengel, C. Shen, J. Gao, F. Hu, and Z. Zhang. Learning graph structure for multi-label image classification via clique generation. In *CVPR*, pages 4100–4109, 2015. 2

G. Tsoumakas, I. Katakis, and I. Vlahavas. Effective and efficient multilabel classification in domains with large number of labels. In *Proc. ECML/PKDD 2008 Workshop on Mining Multidimensional Data*, 2008. 12

G. Tsoumakas, I. Katakis, and I. Vlahavas. *Mining Multi-label Data*, pages 667–685. Springer US, Boston, MA, 2010. 1

V. N. Vapnik. *Statistical learning theory*. Wiley, Wiley New York, 1998. 5, 8

C. Xu, T. Liu, D. Tao, and C. Xu. Local Rademacher complexity for multi-label learning. *IEEE Trans. Image Process.*, 25(3):1495–1507, 2016. 7, 10, 11

H.-F. Yu, P. Jain, P. Kar, and I. S. Dhillon. Large-scale multi-label learning with missing labels. In *ICML*, pages 593–601, 2014. 2, 5, 10, 11, 20

T. Zhang. Statistical behavior and consistency of classification methods based on convex risk minimization. *Ann. Statist.*, 32(1):56–85, 02 2004. 2, 5, 8, 9

# A Appendix

## A.1 Generalization bound

**Theorem 4** (McDiarmid's inequality [McDiarmid, 1989])**.** *Let $Z_1,\dots,Z_N$ be independent random variables taking values in $\mathcal{Z}$. Further, let $\varphi : \mathcal{Z}^N \mapsto \mathbb{R}$ be a function of $Z_1,\dots,Z_N$ that satisfies $\forall i \in [N] \quad \forall z_1,\dots,z_N, z_i' \in \mathcal{Z}$,*

$$\left|\varphi(z_1,\dots,z_i,\dots,z_N) - \varphi(z_1,\dots,z_i',\dots,z_N)\right| \leq c_i \ .$$

*Then for all $t>0$,*

$$\mathbb{P}\{\phi \leq \mathbb{E}\phi + t\} \leq \exp\left(-\frac{2t^2}{\sum_{i=1}^N c_i^2}\right) \ .$$

*Proof of Proposition 1.* Let $\varphi(\mathcal{D}_N) = \sup_{f\in\mathcal{F}}\left|\mathcal{R}(f) - \hat{\mathcal{R}}_N(f,\mathcal{D}_N)\right|$ and denote by $\bar{f} \in \mathcal{F}$ the classifier[3] such that $\varphi(\mathcal{D}_N) = \left|\mathcal{R}(\bar{f}) - \hat{\mathcal{R}}_N(\bar{f},\mathcal{D}_N)\right|$. Let $\mathcal{D}_N', \mathcal{D}_N$ be two sets of observations such that $\mathcal{D}_N'$ and $\mathcal{D}_N$ differ only in one observation $(X_k, Y_k)$ and $(X_k', Y_k')$, where $(X_k', Y_k')$ is a copy of $(X, Y)$ independent from $\mathcal{D}_N$, hence

$$\begin{aligned}
\left|\varphi(\mathcal{D}_N) - \varphi(\mathcal{D}_N')\right| \leq& \frac{p_1}{N}\left|\sum_{l=1}^{L}\left(\mathcal{L}_1(\bar{f}^l(X_k))\mathbb{1}_{\{Y_k^l=1\}} - \mathcal{L}_1(\bar{f}^l(X_k'))\mathbb{1}_{\{(Y')_k^l=1\}}\right)\right| \\
&+\frac{p_0}{N}\left|\sum_{l=1}^{L}\left(\mathcal{L}_0(\bar{f}^l(X_k))\mathbb{1}_{\{Y_k^l=0\}} - \mathcal{L}_0(\bar{f}^l(X_k'))\mathbb{1}_{\{(Y')_k^l=0\}}\right)\right| \\
\leq& \frac{p_1}{N}\left|\sum_{l=1}^{L}\left(\mathcal{L}_1(f^l(X_i))\mathbb{1}_{\{Y_k^l=1\}} - \mathcal{L}_1(\bar{f}^l(X_k'))\mathbb{1}_{\{(Y')_k^l=1\}}\right)\right| + \frac{L}{N}p_0 \\
\leq& \frac{2K}{N}p_1 + \frac{L}{N}p_0 \ ,
\end{aligned}$$

where in the last inequality we used that $\sum_{l=1}^L(\mathcal{L}_1(\bar{f}^l(X_k))\mathbb{1}_{\{Y_k^l=1\}} - \mathcal{L}_1(\bar{f}^l(X_k'))\mathbb{1}_{\{(Y')_k^l=1\}})$ consists of at most $2K$ non-zero terms almost surely under Assumption 1. □

### A.1.1 Refined bound

The following theorem [Bousquet, 2002] is used to prove Lemma 2.

**Theorem 5** ([Bousquet, 2002])**.** *Let $X_1,\dots,X_N$ be a sequence of independent random variables with values in some polish space $\mathcal{X}$ and distributed according to $\mathbb{P}$. Let $F$ be a $\mathbb{P}$-measurable function from $\mathcal{X}^N$ to $\mathbb{R}$. Let $Z = F(X_1,\dots,X_N)$, let $\mathcal{A}_k$ be the sigma algebra generated by $(X_1,\dots,X_k)$ for all $k\in[N]$ and let $\mathcal{A}_n^k$ be the sigma field generated by $(X_1,\dots,X_{k-1},X_{k+1},\dots,X_N)$ for all $k\in[N]$. We denote by $\mathbb{E}_n^k[\cdot]$ the expectation taken conditionally to $\mathcal{A}_n^k$. Let $(Z, Z_1',\dots,Z_N')$, be a sequence of $\mathcal{A}$-measurable random variables and let $(Z_k)_{k=1,\dots,N}$ be a sequence of random variables respectively $\mathcal{A}_n^k$-measurable. Assume that the following inequalities are satisfied*

$$Z_k' \leq Z - Z_k \leq 1 \text{ a.s.}, \quad \mathbb{E}_n^k[Z_k'] \geq 0 \text{ and } Z_k' \leq 1 \text{ a.s.} \tag{24}$$

*Let $r>0$ be a real satisfying $\sum_{k=1}^N \mathbb{E}_n^k[(Z_k')^2] \leq Nr$ almost surely. If the following condition holds*

$$\sum_{k=1}^{N}(Z - Z_k) \leq Z \text{ a.s.}, \tag{25}$$

[3] We suppose that the supremum is reached in the definition of $\varphi(\mathcal{D}_N)$. When this is not the case an $\epsilon$ argument would be needed that we leave to the reader.

*hence for all* $\delta > 0$ *we obtain*

$$\mathbb{P}\left\{Z \le \mathbb{E}[Z] + \sqrt{(4\mathbb{E}[Z]+Nr)\log(1/\delta)} + \frac{\log(1/\delta)}{3}\right\} \le 1-\delta \ . \tag{26}$$

The following result is the consequence of Theorem 5.

**Theorem 6.** *Let* $r > 0$ *be such that*

$$\forall j \in \{0,1\}, \quad \forall f \in \mathcal{F}, \quad \mathbb{E}[\max_{l\in[L]} \mathcal{L}_j^2(f^l(X))] \le r \ , \tag{27}$$

*and*

$$Z = \sup_{f\in\mathcal{F}} \left|\mathcal{R}(f) - \hat{\mathcal{R}}_N(f)\right| \ .$$

*With* $p_0, p_1$ *as in Eq. (8) hence under Assumption (1) with probability at least* $1-\delta$ *we have*

$$Z \le \mathbb{E}[Z] + \sqrt{\mathbb{E}[Z]\frac{8K}{N}\log(1/\delta) + \frac{32K^2r}{N}\log(1/\delta)} + \frac{4K}{3N}\log(1/\delta) \ .$$

*Proof of Theorem 6.* To prove Theorem 6 we need to justify the conditions of Theorem 5, with appropriate choice of random variables. We define the following variables

$$Z = \frac{N}{4K}\sup_{f\in\mathcal{F}} \left|\mathcal{R}(f) - \hat{\mathcal{R}}_N(f,\mathcal{D}_N)\right| \ , \tag{28}$$

and

$$Z_k = \frac{N}{4K}\sup_{f\in\mathcal{F}} \left|\mathcal{R}(f) - \hat{\mathcal{R}}_N(f,\mathcal{D}_N^{-k})\right|, \quad Z_k' = \frac{N}{4K}\left|\mathcal{R}(f_k) - \hat{\mathcal{R}}_N(f_k,\mathcal{D}_N)\right| - Z_k \ , \tag{29}$$

where $\mathcal{D}_N^{-k} = \mathcal{D}_N \setminus \{(X_k, Y_k)\}$ and $f_k$ is such that $Z_k = \frac{N}{4K}\left|\mathcal{R}(f_k) - \hat{\mathcal{R}}_N(f_k,\mathcal{D}_N^{-k})\right|$. We denote by $\bar{f}$ the function for which the supremum is reached in $Z$. Hence we get

$$\begin{aligned}
Z_k' \le Z - Z_k &\le \frac{N}{4K}\left|\mathcal{R}(\bar{f}) - \hat{\mathcal{R}}_N(\bar{f},\mathcal{D}_N)\right| - \frac{N}{4K}\left|\mathcal{R}(\bar{f}) - \hat{\mathcal{R}}_N(\bar{f},\mathcal{D}_N^{-k})\right| \\
&\le \frac{1}{4K}\left|\sum_{l=1}^{L} p_0\mathcal{L}_0(\bar{f}^l(X_k))\mathbb{1}_{\{Y_k^l=0\}} + p_1\mathcal{L}_1(\bar{f}^l(X_k))\mathbb{1}_{\{Y_k^l=1\}}\right| \\
&\le 1 \quad \text{a.s.} \ ,
\end{aligned}$$

where in the last inequality we used Assumption 1 and boundedness of $\mathcal{L}_0$ and $\mathcal{L}_1$. Moreover we have

$$\mathbb{E}_n^k[Z_k'] \ge \frac{N}{4K}\left|\mathbb{E}_n^k\left[\mathcal{R}(f_k) - \hat{\mathcal{R}}_N(f_k,\mathcal{D}_N)\right]\right| - Z_k = 0 \ ,$$

We also have

$$(N-1)Z = \frac{N}{4K}\left|\sum_{k=1}^{N}(\mathcal{R}(\bar{f}) - \hat{\mathcal{R}}_N(\bar{f},\mathcal{D}_N^{-k}))\right| \le \frac{N}{4K}\sum_{k=1}^{N}\left|\mathcal{R}(\bar{f}) - \hat{\mathcal{R}}_N(\bar{f},\mathcal{D}_N^{-k})\right| \le \sum_{k=1}^{N} Z_k \ ,$$

hence

$$\sum_{k=1}^{N}(Z - Z_k) \le Z \ ,$$

and finally we have

$$\begin{aligned}\sum_{k=1}^{N}\mathbb{E}_n^k[(Z_k')^2] &\leq \sum_{k=1}^{N}\mathbb{E}_n^k\Big[\Big(\frac{1}{4K}\sum_{l=1}^{L}p_0\mathcal{L}_0(f_k^l(X_k))\mathbb{1}_{\{Y_k^l=0\}}+p_1\mathcal{L}_1(f_k^l(X_k))\mathbb{1}_{\{Y_i^l=1\}}\Big)^2\Big]\\ &\leq \frac{1}{16K^2}\sum_{k=1}^{N}\mathbb{E}_n^k\Big[\Big(p_0L\max_{l\in[L]}\mathcal{L}_0(f_k^l(X_k))+Kp_1\max_{l\in[L]}\mathcal{L}_1(f_k^l(X_k))\Big)^2\Big]\\ &\leq \frac{1}{8K^2}\sum_{k=1}^{N}\mathbb{E}_n^k\Big[p_0^2L^2\max_{l\in[L]}\mathcal{L}_0^2(f_k^l(X_k))\Big]+\mathbb{E}_n^k\Big[K^2p_1^2\max_{l\in[L]}\mathcal{L}_1^2(f_k^l(X_k))\Big]\\ &\leq \frac{Nr}{8K^2}\left(p_0^2L^2+K^2p_1^2\right)=\frac{Nr}{8}\left(p_0^2\frac{L^2}{K^2}+p_1^2\right)\leq Nr \quad \text{a.s.}\ ,\end{aligned}$$

and we conclude. □

*Proof of Lemma 2.* It is sufficient to apply the following elementary result to the Theorem 6, to obtain the proof of Lemma 2.

**Lemma 4.** For every $u, v \geq 0$,

$$\sqrt{u+v} \leq \sqrt{u}+\sqrt{v}, \tag{30}$$

$$2\sqrt{uv} \leq u+v. \tag{31}$$

Using previous inequality to bound product we can rewrite Equation (28) as follows

$$Z \leq 2\mathbb{E}[Z]+K\sqrt{\frac{32r}{N}\log(1/\delta)}+K\frac{10}{3N}\log(1/\delta)\ . \tag{32}$$

□

### A.2 On convexification

*Proof of Lemma 3.*

$$\begin{aligned}\mathcal{R}^w(f)-\mathcal{R}^{w*} &= \mathbb{E}_X\sum_{l=1}^{L}p_0(1-\eta^l(X))(\mathbb{1}_{\{f^l(X)\geq 0\}}-\mathbb{1}_{\{\eta^l(X)\geq p_0\}})+p_1\eta^l(X)(\underbrace{\mathbb{1}_{\{f^l(X)<0\}}-\mathbb{1}_{\{\eta^l(X)<p_0\}}}_{=\mathbb{1}_{\{\eta^l(X)\geq p_0\}}-\mathbb{1}_{\{f^l(X)\geq 0\}}})\\ &= \sum_{l=1}^{L}\mathbb{E}_X(p_0-\eta^l(X))(\mathbb{1}_{\{f^l(X)\geq 0\}}-\mathbb{1}_{\{\eta^l(X)\geq p_0\}})\\ &\leq \sum_{l=1}^{L}\mathbb{E}_X[\mathbb{1}_{\{(\eta^l(X)-p_0)f^l(X)\leq 0\}}\left|\eta^l(X)-p_0\right|]\\ &\leq \sum_{l=1}^{L}\left(\mathbb{E}_X[\mathbb{1}_{\{(\eta^l(X)-p_0)f^l(X)\leq 0\}}\left|\eta^l(X)-p_0\right|^s]\right)^{1/s} \qquad \text{(by Jensen's inequality)}\\ &\leq c\sum_{l=1}^{L}\left(\mathbb{E}_X[\mathbb{1}_{\{(\eta^l(X)-p_0)f^l(X)\leq 0\}}\Delta\mathcal{R}_\phi^w(\eta^l,0)]\right)^{1/s}\ . \qquad \text{(thanks to Eq. (20))}\end{aligned}$$

Now it is sufficient to show that for all $l \in [L]$ having $(\eta^l(X)-p_0)f^l(X) \leq 0$ implies $\Delta\mathcal{R}_\phi^w(\eta^l,0) \leq \Delta\mathcal{R}_\phi^w(\eta^l,f^l)$. We only need to prove that $(\eta^l-p_0)p^l \leq 0$ implies $\mathcal{R}_\phi^w(\eta^l,0) \leq \mathcal{R}_\phi^w(\eta^l,p^l)$. To see this, we consider the following cases:

- $\eta > p_0$: By assumption, we have $(f^*_\phi)^l(\eta^l) > 0$. Moreover, $(\eta - p_0)p^l \leq 0$ implies $p^l \leq 0$. Hence $0 \in [p^l, (f^*_\phi)^l(\eta^l)]$ and by convexity of the real function $t \mapsto \mathcal{R}^w_\phi(\eta^l, t)$ we have $\mathcal{R}^w_\phi(\eta^l, 0) \leq \max\{\mathcal{R}^w_\phi(\eta^l, p^l), \mathcal{R}^w_\phi(\eta^l, (f^*_\phi)^l(\eta^l))\} = \mathcal{R}^w_\phi(\eta^l, p^l)$.
- $\eta < p_0$: Note that we require $(f^*_\phi)^l(\eta^l) < 0$. Therefore the same analysis as above brings the desired result.
- $\eta^l = p_0$: Since $(f^*_\phi)^l(\eta^l) = 0$ we have the desired result.

Hence we have proved the following inequality

$$\mathcal{R}^w(f) - \mathcal{R}^{w*} \leq c \sum_{l=1}^{L} \Big(\mathbb{E}_X \Delta\mathcal{R}^w_\phi(\eta^l, f^l)\Big)^{1/s} \ .$$

Then, by concavity of $x^{1/s}$ we conclude that

$$\mathcal{R}^w(f) - \mathcal{R}^{w*} \leq cL^{1-1/s}\Big(\sum_{l=1}^{L} \mathbb{E}_X \Delta\mathcal{R}^w_\phi(\eta^l, f^l)\Big)^{1/s} = cL^{1-1/s}\Big(\mathcal{R}^w_\phi(f) - \mathcal{R}^w_\phi(f^*_\phi)\Big)^{1/s} \ .$$

□

*Proof for Examples 1.* We first notice, that the minimization in Equation (19) is separable, *i.e.,* $f^*_\phi(\eta) = ((f^*_\phi)^1(\eta^1), \ldots, (f^*_\phi)^L(\eta^L))^\top$, where for all $l \in [L]$

$$(f^*_\phi)^l(\eta^l) \in \arg\min_f \mathcal{R}^w_\phi(\eta^l, f^l) \ . \tag{33}$$

- Let $\phi(v) = (1-v)^2$ and fix $l \in [L]$, minimizing $\mathcal{R}^w_\phi(\eta^l, f^l)$ over $f^l$ we obtain

  $$(f^*_\phi)^l(\eta^l) = \frac{p_0 - \eta^l}{\eta^l(1-2p_0) + p_0} \ ,$$

  substituting to $\Delta\mathcal{R}^w_\phi(\eta^l, 0)$ we get

  $$\Delta\mathcal{R}^w_\phi(\eta^l, 0) = \frac{\left|p_0 - \eta^l\right|^2}{\eta^l(1-2p_0) + p_0} \ ,$$

  and finally, since $\eta^l(1-2p_0) + p_0 \leq 2$ we conclude by

  $$\left|p_0 - \eta^l\right|^2 \leq 2\frac{\left|p_0 - \eta^l\right|^2}{\eta^l(1-2p_0) + p_0} \ .$$

- Let $\phi(v) = e^{-v}$ and fix $l \in [L]$, minimizing $\mathcal{R}^w_\phi(\eta^l, f^l)$ over $f^l$ we obtain

  $$(f^*_\phi)^l(\eta^l) = \frac{1}{2}\log\Big(\frac{\eta^l(1-p_0)}{p_0(1-\eta^l)}\Big) \ ,$$

  substituting to $\Delta\mathcal{R}^w_\phi(\eta^l, 0)$ we get

  $$\Delta\mathcal{R}^w_\phi(\eta^l, 0) = \left|\sqrt{(1-p_0)\eta^l} - \sqrt{p_0(1-\eta^l)}\right|^2 \ ,$$

and finally, since $\left|\sqrt{(1-p_0)\eta^l}+\sqrt{p_0(1-\eta^l)}\right|^2 \leq 2$ we conclude by

$$\begin{aligned}\left|p_0-\eta^l\right|^2 &= \left|\sqrt{(1-p_0)\eta^l}-\sqrt{p_0(1-\eta^l)}\right|^2\left|\sqrt{(1-p_0)\eta^l}+\sqrt{p_0(1-\eta^l)}\right|^2 \\ &\leq 2\left|\sqrt{(1-p_0)\eta^l}-\sqrt{p_0(1-\eta^l)}\right|^2 \ .\end{aligned}$$

□

### A.3 Linear prediction approach

*Proof of Theorem 3.* Our proof follows the analysis presented in [Yu et al., 2014] with slight modification and use of Assumption 1. It is sufficient to provide bound for

$$\mathfrak{R}_1(\mathcal{W}) = \frac{1}{N}\mathbb{E}_N\mathbb{E}_\varepsilon \sup_{f\in\mathcal{F}}\left\{\sum_{i=1}^N\sum_{l=1}^L \varepsilon_i^l f^l(X_i)\mathbb{1}_{\{Y_i^l=1\}}\right\} \ ,$$

$$\mathfrak{R}_0(\mathcal{W}) = \frac{1}{N}\mathbb{E}_N\mathbb{E}_\varepsilon \sup_{f\in\mathcal{F}}\left\{\sum_{i=1}^N\sum_{l=1}^L \varepsilon_i^l f^l(X_i)\mathbb{1}_{\{Y_i^l=0\}}\right\} \ ,$$

to this end let $\boldsymbol{V}_1=[V_1^1,\ldots,V_1^L]$, where $V_1^l=\sum_{i=1}^N \varepsilon_i^l X_i\mathbb{1}_{\{Y_i^l=1\}}$, hence

$$\begin{aligned}\mathfrak{R}_1(\mathcal{W}) &= \frac{1}{N}\mathbb{E}_N\mathbb{E}_\varepsilon \sup_{\boldsymbol{W}\in\mathcal{W}}\langle \boldsymbol{W},\boldsymbol{V}_1\rangle \leq \frac{\lambda}{N}\mathbb{E}_N\mathbb{E}_\varepsilon\|\boldsymbol{V}_1\|_{\sigma,\infty} \\ &\leq \frac{\lambda}{N}\sqrt{\mathbb{E}_N\mathbb{E}_\varepsilon\|\boldsymbol{V}_1\|_{\sigma,2}^2} = \frac{\lambda}{N}\sqrt{\mathbb{E}_N\mathbb{E}_\varepsilon\sum_{l=1}^L\left\|V_1^l\right\|_2^2} \ ,\end{aligned}$$

using expression for $V_1^l$ we can write

$$\begin{aligned}\mathbb{E}_N\mathbb{E}_\varepsilon\sum_{l=1}^L\left\|\sum_{i=1}^N\varepsilon_i^l X_i\mathbb{1}_{\{Y_i^l=1\}}\right\|_2^2 &= \mathbb{E}_N\mathbb{E}_\varepsilon\sum_{l=1}^L\sum_{i=1}^N\|X_i\|_2^2\mathbb{1}_{\{Y_i^l=1\}}+\mathbb{E}_N\mathbb{E}_\varepsilon\sum_{l=1}^L\sum_{i\neq i'}\varepsilon_i^l\varepsilon_{i'}^l\mathbb{1}_{\{Y_i^l=1\}}\mathbb{1}_{\{Y_{i'}^l=1\}}\langle X_i,X_{i'}\rangle \\ &= \mathbb{E}_N\sum_{i=1}^N\|X_i\|_2^2\sum_{l=1}^L\mathbb{1}_{\{Y_i^l=1\}} \leq K\mathbb{E}_N\sum_{i=1}^N\|X_i\|_2^2 \leq KN \ .\end{aligned}$$

where in the second equality we used the definition of Rademacher variables and in the second last inequality we used the fact that $\sum_{l=1}^L\mathbb{1}_{\{Y_i^l=1\}}\leq K$ almost surely under Assumption 1, therefore we have

$$\mathfrak{R}_1(\mathcal{W}) \leq \lambda\sqrt{\frac{K}{N}} \ .$$

For the second Rademacher complexity we can write

$$\begin{aligned}\mathfrak{R}_0(\mathcal{W}) &= \frac{1}{N}\mathbb{E}_N\mathbb{E}_\varepsilon \sup_{\boldsymbol{W}\in\mathcal{W}}\langle \boldsymbol{W},\boldsymbol{V}_0\rangle \leq \frac{\lambda}{N}\mathbb{E}_N\mathbb{E}_\varepsilon\|\boldsymbol{V}_0\|_{\sigma,\infty} \\ &\leq \frac{\lambda}{N}\sqrt{\mathbb{E}_N\mathbb{E}_\varepsilon\|\boldsymbol{V}_0\|_{\sigma,2}^2} = \frac{\lambda}{N}\sqrt{\mathbb{E}_N\mathbb{E}_\varepsilon\sum_{l=1}^L\left\|V_0^l\right\|_2^2} \ ,\end{aligned}$$

**Algorithm 1** Frank-Wolfe algorithm for minimizing $F(X)$ with a nuclear norm constraint.

**Result**: $\hat{\boldsymbol{W}}$
$k = 0$
$\boldsymbol{W}_0 = 0$
**while** *Stopping criteria is not achieved* **do**
    $[u_k, s_k, v_k] = -\mathrm{svds}(-\nabla F(\boldsymbol{W}_k), 1)$
    $\alpha_k = \frac{2}{k+2}$
    $\boldsymbol{W}_{k+1} = (1-\alpha_k)\cdot \boldsymbol{W}_k + \alpha_k \cdot \lambda \cdot u_k v_k^\top$
    $k = k+1$

where $\boldsymbol{V}_0 = [V_0^1, \dots, V_0^L]$, where $V_0^l = \sum_{i=1}^N \varepsilon_i^l X_i \mathbb{1}_{\{Y_i^l = 0\}}$, hence

$$\mathbb{E}_N \mathbb{E}_\varepsilon \sum_{l=1}^L \left\| \sum_{i=1}^N \varepsilon_i^l X_i \right\|_2^2 = \mathbb{E}_N \mathbb{E}_\varepsilon \sum_{l=1}^L \sum_{i=1}^N \|X_i\|_2^2 \mathbb{1}_{\{Y_i^l=0\}} + \mathbb{E}_N \mathbb{E}_\varepsilon \sum_{i \neq i'} \varepsilon_i^l \varepsilon_{i'}^l \mathbb{1}_{\{Y_i^l=0\}} \mathbb{1}_{\{Y_{i'}^l=0\}} \langle X_i, X_{i'} \rangle$$
$$= \mathbb{E}_N \sum_{i=1}^N \|X_i\|_2^2 \sum_{l=1}^L \mathbb{1}_{\{Y_i^l=0\}} \leq L \mathbb{E}_N \sum_{i=1}^N \|X_i\|_2^2 \leq LN \;,$$

Hence we have

$$\mathfrak{R}_1(\mathcal{W}) \leq \lambda \sqrt{\frac{K}{N}} \;, \tag{34}$$

$$\mathfrak{R}_0(\mathcal{W}) \leq \lambda \sqrt{\frac{L}{N}} \;. \tag{35}$$

We conclude by substituting $p_0, p_1$ and Equation (34) into Theorem 1 and we get with probability at least $1-\delta$

$$\sup_{\boldsymbol{W} \in \mathcal{W}} \left\{ \mathcal{R}(\boldsymbol{W}) - \hat{\mathcal{R}}_N(\boldsymbol{W}) \right\} \leq 2C\lambda \left( \sqrt{\frac{K}{N}} + \frac{K}{\sqrt{NL}} \right) + 4\frac{K}{\sqrt{N}} \sqrt{\log 1/\delta} \;. \tag{36}$$

□

## A.4 Frank-Wolfe algorithm

Let $F : \mathbb{R}^{D\times L} \mapsto \mathbb{R}$ be convex function, hence we would like to solve the following minimization problem

$$\min_{\boldsymbol{W} \in \mathbb{R}^{D\times L}} \quad F(\boldsymbol{W})$$
$$\text{s.t.} \quad \|\boldsymbol{W}\|_{\sigma,1} \leq \lambda.$$

One of the possible approach is to consider the Frank-Wolfe[4] algorithm [Frank and Wolfe, 1956] (see also [Jaggi, 2013, Jaggi and Sulovský, 2010] for applications to similar matrix optimization problems). This algorithm is attractive due to its cheap iteration cost, *i.e.,* on each iteration we need to evaluate only the top singular value of the gradient of $F$ (in contrast to forward-backward variants that would require a full SVD evaluation at each iteration). The Frank-Wolfe algorithm applied to Problem (A.4) is summarized in Algorithm 1. The function $\mathrm{svds}(-\nabla F(\boldsymbol{W}_k), 1)$ computes the largest singular value and the corresponding singular vectors.

[4] also referred to as conditional gradient

As a stopping criterion we follow the one proposed in [Jaggi, 2013] that can be specified in our context as:

$$|\langle \boldsymbol{W}_k, \nabla F(\boldsymbol{W}_k)\rangle + \lambda \cdot s_k| \leq \epsilon, \tag{37}$$

for some $\epsilon > 0$ This duality gap criterion ensures that $F(\boldsymbol{W}_k) - F(\hat{\boldsymbol{W}}) \leq \epsilon$ for the output of the algorithm.

### A.5 Synthetic data

First, we remind that the sigmoid function is given by:

$$S(t) = \frac{1}{1+e^{-t}} \quad . \tag{38}$$

We consider the following framework to simulate synthetic datasets with controlled sparsity constant $K$. Each entry of the design matrix $\boldsymbol{X} \in \mathbb{R}^{N \times D}$ is simulated according to a uniform distribution on $[-1, 1]$. We form a matrix of coefficients $\boldsymbol{B} \in \mathbb{R}^{L \times D}$, in which every entry is either 2 or $-2$ with probability 0.5. We transform the matrix $\boldsymbol{Z} = \boldsymbol{X}\boldsymbol{B}^\top$ into probabilities by applying the sigmoid function Eq. (38) to each entry. Having the matrix $\boldsymbol{Z}$ of probabilities we pick the $K/2$ highest values in each row of $\boldsymbol{Z}$ and set them to 1. Moreover, we set the $L-K$ lowest values in each row of $\boldsymbol{Z}$ to 0. Now, we have a matrix with $K/2$ ones, $K/2$ value between zero and one (probabilities) and $L-K$ zeros in each row. Each element of the label matrix $Y$ is obtained as Bernoulli random variable with probability from the described matrix, see Equation (39).

$$Y_i^l = \begin{cases} 1, & \text{with probability } Z_i^l \\ 0, & \text{with probability } 1 - Z_i^l \end{cases} \quad . \tag{39}$$

This framework, allows to obtain $\boldsymbol{Y}$, in which rows are at least $K/2$ sparse and at most $K$ sparse.